\documentclass[12pt,leqno]{amsart}
\usepackage{amsmath}
\usepackage{amssymb}

 \newtheorem{thm}{Theorem}[section]
 
 \newtheorem{lem}[thm]{Lemma}
 
 \theoremstyle{definition}
 \newtheorem{defn}[thm]{Definition}
 \theoremstyle{remark}

 \numberwithin{equation}{section}

\begin{document}

\title{On A Fragment of the Universal Baire Property for $\Sigma^1_2$ sets}

\author{Stuart Zoble \\ March 18, 2006}

\address{Department of Mathematics, University of Toronto}
\email{szoble@math.toronto.edu}



\begin{abstract}  There is a well-known global equivalence between $\Sigma^1_2$ sets having 
the Universal Baire property, two-step $\Sigma^1_3$ generic absoluteness, and the closure 
of the universe under the sharp operation.  In this note, we determine the exact consistency strength 
of $\Sigma^1_2$ sets being $(2^{\omega})^{+}$-cc Universally Baire, which is below $0^{\#}$.  
In a model obtained, there is a $\Sigma^1_2$ set which is weakly $\omega_2$-Universally Baire but not 
$\omega_2$-Universally Baire.  
\end{abstract}

\maketitle

\section{Introduction} \noindent Consider the following two properties of a set of reals 
$A \subset \omega^{\omega}$ at an infinite cardinal $\kappa$.

\begin{enumerate}

\item For every continuous $f:\kappa^{\omega} \rightarrow \omega^{\omega}$ there is a dense set of 
$p \in \kappa^{<\omega}$ such that either $f^{-1}(A)$ is meager or comeager below $p$. 

\item For every continuous $f:\kappa^{\omega} \rightarrow \omega^{\omega}$ there is a dense set of 
$p \in \kappa^{<\omega}$ such that either $f^{-1}(A) \cap \sigma^{\omega}$ is meager below $p$ in 
$\sigma^{\omega}$ for a club of $\sigma \in [\kappa]^{\omega}$ or comeager below $p$ in 
$\sigma^{\omega}$ for a club of $\sigma \in [\kappa]^{\omega}$.

\end{enumerate}

\noindent The first property asserts that $A$ is $\kappa$-Universally Baire or fully captured at $\kappa$, 
and the second that $A$ is weakly $\kappa$-Universally Baire (the author's coinage) or weakly captured at 
$\kappa$ (see [6] and Lemma 4.1 below).  
The implication (1) to (2) is immediate as a club of $\sigma \in [\kappa]^{\omega}$ are closed under a 
Banach-Mazur strategy in the space $\kappa^{\omega}$.  Regarding the reverse implication, any set 
of reals of size $\omega_1$ is a counterexample at $\kappa = \omega_2$ assuming Martin's Maximum 
(see Thm. 2.6 of [6] and Thm. 3.1 of [4]).  The place to look for a definable counterexample is the 
pointclass $\Sigma^1_2$ with $\kappa$ either $\omega_1$ or $\omega_2$.  This is because 
(1) and (2) are equivalent for $\kappa = \omega$ and for $\Delta^1_2$ sets as a whole.  As the 
particular scenario for a counterexample suggested by [6] involves the question of the 
consistency strength of $\Sigma^1_2$ sets being $\omega_1$ or $\omega_2$-Universally Baire, in 
particular whether this is possible without sharps, this paper was motivated by the following 
question.  What fragment of the Universal Baire Property can $\Sigma^1_2$ sets
have below $0^{\#}$, as measured by the weight or cellularity of the
preimage space?  In Theorem 3.4 of [1] a global equivalence is established between
$\Sigma^{1}_{2}$ sets being Universally Baire, two-step $\Sigma^1_3$
Generic Absoluteness, and the closure of the universe under the
sharp operation. This equivalence, however, is not true
level-by-level.  In particular, the relevant part of their argument (originating in [4]) would
require that $\Sigma^{1}_{2}$ sets be $\omega_{\omega+1}$-Universally
Baire to prove the existence of $0^{\#}$.  Using the full strength
of Covering for $L$ this can be reduced to $\omega_3$.  On the other hand, 
Woodin has shown that $\Pi^1_2$ sets can be $\omega_2$-cc-Universally Baire in a forcing 
extension of $L$.    

\begin{defn} $A \subset \mathbb{R}$ is $\kappa$-cc-Universally Baire if $f^{-1}(A)$ 
has the Baire property in $X$ for every completely regular space $X$ of cellularity less than 
$\kappa$ and continuous map $f:X \rightarrow \mathbb{R}$.
\end{defn}

\begin{thm} (Woodin) Assume $\lambda_0 < \lambda_1 < \lambda_2$ are cardinals of $L$ 
and there is an elementary embedding $\pi:L_{\lambda_1} \rightarrow L_{\lambda_2}$ 
with critical point $\lambda_{0}$ such that $\pi(\lambda_{0}) = \lambda_1$.  Then $\Sigma^1_2$ 
sets are $\omega_2$-cc-Universally Baire in a forcing extension of $L$ in which $CH$ holds.
\end{thm}

\noindent For equivalent versions of Definition 1.1, the reader is referred to Theorem 2.1 of [1].  
In particular, we will use that a $\kappa$-cc Universally Baire set remains ccc-Universally Baire 
after forcing with a $\kappa$-cc poset.  When combined with the argument of Theorem 3.4 of [1], the 
above shows that two-step $\Sigma^1_3$ generic absoluteness for $(2^{\omega})^{+}$-cc forcings 
$\mathbb{P} * \dot{\mathbb{Q}}$ can hold in a forcing extension of $L$.  In this note we 
reduce the hypothesis of Theorem 1.2 to obtain an exact equiconsistency.  

\begin{defn} An ordinal $\kappa$ is $L$-large to $\lambda$ if for
every $\alpha < \lambda$ there is an elementary $j:L_{\alpha}
\rightarrow L_{\beta}$ with critical point $\kappa$ such that
$j(\kappa) \geq \alpha$.  \end{defn}

\noindent Note that an ordinal $\kappa$ is $L$-large to some
cardinal of $L$ if and only if $\kappa$ is weakly compact in $L$. We
obtain stronger notions by requiring that $\lambda$ be inaccessible or 
weakly compact in L.  Neither notion implies that $0^{\#}$
exists (simply collapse $\lambda^{+}$ and use absoluteness of $L[g]$), though $L$ cannot see such an 
embedding with $\alpha \geq (\kappa^{+})^{L}$.

\begin{thm} The following are equiconsistent.
\begin{enumerate}
\item[(a)] $\Sigma^1_2$ sets are $(2^{\omega})^{+}$-cc-Universally Baire
\item[(b)] There is a $\kappa$ which is $L$-large to a weakly
compact of $L$.
\end{enumerate}
\end{thm}
\noindent In the last section of 
this paper, we argue that in the model of Theorem 1.3 (a) there must be a $\Sigma^1_2$ set which has 
a weak capturing term at $\omega_2$ but no full capturing term at $\omega_2$.

\begin{thm}  It is consistent that $\Sigma^1_2$ sets are weakly $\omega_2$-Universally Baire but not 
$\omega_2$-Universally Baire.
\end{thm}

\noindent In what follows, any space $X^{\omega}$ will carry the product topology, with the set $X$ 
given the discrete topology.  All pointclasses are boldface, and every statement below involving 
$\Sigma^1_2$ sets applies equally well to $\Pi^1_2$ sets.  We would like to thank Hugh Woodin 
for sharing his proof of Theorem 1.2 and allowing us to include elements of it here, and Stevo 
Todorcevic for several helpful comments regarding earlier drafts of this paper.
 
\section{Preliminaries}

\noindent Call a term $\dot{A}$ a $Col(\omega,\kappa)$ capturing term for a
set of reals $A$ if there is a club of countable elementary
submodels $X \prec H(\theta)$ with transitivization $H$ and collapse
map $\pi$ such that$$\pi(\dot{A})_{g} = A \cap H[g]$$

\noindent for every $H$-generic $g \subset Col(\omega,
\pi(\kappa))$.  A set of reals $A$ has a $Col(\omega,\kappa)$
capturing term if and only if $A$ is $\kappa$-Universally Baire (see
Lemma 1.6 of [6]). The proof below uses this observation and an
argument from [1].

\begin{thm}  The following are equivalent for a cardinal $\kappa$.

\begin{enumerate}
\item $\Sigma^{1}_{2}$ sets are $\kappa$-Universally Baire

\item For all sufficiently large $\theta$, there is a club of
countable $X \prec H(\theta)$ such that $X[g]$ is $\Sigma^1_2$
elementary in $V$ for every $X$-generic $g \subset Col(\omega,\kappa
\cap X)$.
\end{enumerate}
\end{thm}
\begin{proof}  Let $A$ be $\Sigma^1_2$ defined by a formula $\phi$ (we suppress any parameter).
Assuming (2) let $\dot{A}$ be the set of pairs $(p,\tau)$ such that
$p \Vdash \phi(\tau)$.  Thus $\dot{A}$ is a capturing term for $A$
and it follows that $A$ is $\kappa$-Universally Baire.  For the
other direction assume $S$ and $T$ are trees for a $\Sigma^1_2$ set
and its complement.  Suppose the $\Sigma^{1}_{2}$ set is defined by
a formula $\phi(x)$ (again suppressing parameters).  The argument of
Theorem 3.4 of [1] shows that $$p[S]^{V[G]} = \{x \ | \ \phi(x)
\}^{V[G]},$$

\noindent where $G \subset Col(\omega,\kappa)$ is $V$-generic.  This
uses $\Pi^1_1$-Uniformization.  Let $\phi$ be a $\Sigma^1_{2}$
formula defining the universal $\Sigma^1_2$ set $A$ and let $X \prec
H(\theta)$ contain $S$ and $T$. Then $p[S]^{X[g]} = A \cap X[g]$ and
$X[g]$ thinks $p[S]^{X[g]}$ is the universal $\Sigma^1_2$ set. Hence
$X[g]$ is $\Sigma^1_2$ elementary in $V$.  \end{proof}

\noindent  Using (2) and the full strength of Covering for $L$ we may
now argue that $\Sigma^1_2$ sets being $\omega_3$-Universally Baire implies that 
$0^{\#}$ exists.

\begin{thm}  Assume $\Sigma^1_2$ sets are $\omega_3$-Universally Baire.
 Then $0^{\#}$ exists.  \end{thm}
\begin{proof}  Let $\kappa = \omega_3$.  Let $\kappa = \omega_3$.  We first argue that there are club many
$\sigma \in [\kappa]^{\omega}$ such that $otp(\sigma)$ is a regular
cardinal of $L$.  Let $\kappa \in X \prec H(\theta)$ be countable
with transitive collapse $\pi: X \rightarrow \bar{X}$.  Note that
$\pi(\kappa) = otp(X \cap \kappa)$ and that there are club many such
$X \cap \kappa$.  Thus $\bar{X}$ thinks that $\pi(\kappa)$ is
cardinal of $L$.  If $\pi(\kappa)$ were not a regular cardinal of
$L$ then there would be a countable $L_{\gamma}$ which sees this.
Let $g \subset Col(\omega,\pi(\kappa))$ be $\bar{X}$-generic, and
let $z \in X[g]$ be a real coding a well-ordering of length
$\pi(\kappa)$. Then by Theorem 0.4 $\bar{X}[g]$ thinks there is a
level of $L$ which sees that the ordinal coded by $z$ is not
regular. This is a contradiction as $L^{\bar{X}} = L^{\bar{X}[g]}$.
We now argue that the set of $\alpha < \kappa$ such that $\alpha$ is
a regular cardinal of $L$ contains a club in $V$.  Let
$f:\kappa^{<\omega} \rightarrow \kappa$ be such that any $\sigma \in
[\kappa]^{\omega}$ which is closed under $f$ has the property that
$otp(\sigma)$ is a regular cardinal of $L$.  Let $\alpha < \kappa$
such that $f[\alpha^{<\omega}] \subset \alpha$.  Since there are
club many such $\alpha$ it suffices to show that $\alpha$ is a
regular cardinal of $L$.  Suppose not.  Then there is a countable $X
\prec H(\kappa)$ with $\alpha \in X$ such that $X \cap \alpha =
\sigma$ is closed under $f$.  Let $\bar{X}$ be the transitivization
of $X$ with collapse map $\pi$.  Then $\bar{X}$ thinks that
$\pi(\alpha)$ is not a regular cardinal of $L$ hence $\pi(\alpha)$
is not a regular cardinal of $L$ by absoluteness. This contradicts
$\pi(\alpha) = otp(\sigma)$.  It follows that there is an $\alpha$
with $cf(\alpha) < \omega_2 < \alpha < \omega_3$ which is a regular
cardinal of $L$. Let $\sigma \subset \alpha$ be unbounded in
$\alpha$ and have size $cf(\alpha)$. Then $\sigma$ cannot be covered
by a set in $L$ of size $\omega_1$. \end{proof}

\noindent We conjecture that $\omega_3$-cc Universally Baire suffices for the 
argument above.  Under this assumption $\omega_3$ is weakly compact in $L$ by 
Lemma 4 of [5] and a theorem in [2].  We close this section with an equivalence between 
$\Sigma^1_2$ sets being
$\omega_1$-UB and the existence of a club of suitably closed
submodels.  We say that $\omega_2$ is inaccessible to $P(\omega_1)$
if $\omega_2$ is an inaccessible cardinal in $L[X]$ for every $X
\subseteq \omega_1$.

\begin{lem} The following are equivalent.

\begin{enumerate}
\item $\omega_2$ is inaccessible to $P(\omega_1)$ and $\Sigma^{1}_{2}$ sets are
$\omega_1$-Universally Baire

\item $\omega_2$ is
inaccessible to $P(\omega_1)$ and for sufficiently large $\theta$
there is a club of $X \prec H(\theta)$ such that for every $\tau \in
P(\omega_1) \cap X$ and every $L[\tau]$-cardinal $\gamma \in X \cap
\omega_2$ the order type of $\gamma \cap X$ is itself an $L[\tau
\cap X]$-cardinal.

\item For sufficiently large $\theta$ there is a club of $X \prec
H(\theta)$ such that for every $\tau \in P(\omega_1) \cap X$ the
order type of $X \cap \omega_2$ is an $L[\tau \cap X]$-cardinal.

\end{enumerate}
\end{lem}

\begin{proof}  By the argument of Theorem 2.2, condition (3) implies
that $\omega_2$ is inaccessible to $P(\omega_1)$.  Thus (2) and (3)
are equivalent.  Again by a boldface version of an argument from
Theorem 2.2, (1) implies (2).  Let $X \prec H(\theta)$ be as in (2).
Let $\pi:X \rightarrow \bar{X}$ be the collapse map and let $g
\subset Col(\omega,\omega_{1} \cap X)$ be $\bar{X}$-generic. Let $y
= \pi(\tau)_{g}$ be a real in $\bar{X}[g]$. Since $g$ is also
$L[\pi(\tau)]$-generic we have
$$\pi(\omega_{2}) > (\pi(\omega_1)^{+})^{L[\pi(\tau)]} \geq
(\omega_1)^{L[y]}$$ \noindent so that $\bar{X}[g]$ is correct about
$\Sigma^1_2$ facts in the parameter $y$.
\end{proof}

\section{Equiconsistency results}

\noindent Fix a surjection $f_{\gamma}:\omega_1 \rightarrow \gamma$ for
each $\gamma$ between $\omega_1$ and $\omega_2$.  If $\gamma$ is a cardinal of $L$
let $S_{\gamma}$ denote the set of $\alpha < \omega_1$ such that the
order type of $f_{\gamma}[\alpha]$ is an $L$-cardinal. Let $S$ be
the set of $\sigma \in [\omega_{2}]^{\omega}$ such that
$$\sigma \cap \omega_{1} \in \bigcap_{\gamma \in \sigma}
S_{\gamma}.$$

\noindent If we assume that $\omega_2$ is inaccessible in $L$ and
that there are stationary many $\sigma \in [\omega_2]^{\omega}$ such
that $otp(\sigma)$ is an $L$-cardinal then it follows that $S$ is
stationary.  Now let $\mathbb{Q}$ be the countable support product
of $\mathbb{Q}_{\gamma}$, ranging over ordinals $\gamma < \omega_2$ which are $L$-cardinals,  
where $\mathbb{Q}_{\gamma}$ is the poset
for shooting a club through $S_{\gamma}$ with countable conditions.
It follows that $\mathbb{Q}$ is $(\omega,\infty)$-distributive.  If
$CH$ holds then $\mathbb{Q}$ satisfies the $\omega_2$-chain
condition.  The following key lemma is implicit in Woodin's proof of Theorem 1.2. 

\begin{lem} Suppose

\begin{enumerate}

\item Every subset of $\omega_1$ is $L$-generic for some poset $\mathbb{P} \in L$
with $|\mathbb{P}| < \omega_2$

\item There are stationary many $\sigma \in [\omega_2]^{\omega}$
such that $otp(\sigma)$ is an $L$-cardinal

\item $\omega_2$ is inaccessible in $L$ and $CH$ holds.

\end{enumerate}

\noindent Then $\Sigma^1_2$ sets are $\omega_1$-Universally Baire in $V[G]$ where $G
\subset \mathbb{Q}$ is $V$-generic.

\end{lem}

\begin{proof}  We show that condition (3) of Lemma 2.3 is
satisfied in $V[G]$.  As discussed above, $\mathbb{Q}$ preserves
cardinals under these hypotheses and by design there is in $V[G]$ a club of $\sigma \in
[\omega_2]^{\omega}$ such that $otp(\sigma)$ is a cardinal of $L$.
Further, condition (1) continues to hold in $V[G]$.  Suppose $X
\prec H(\theta)$ is such that $otp(X \cap \omega_2)$ is an
$L$-cardinal.  Let $\tau \in P(\omega_1) \cap X$.  Then there are
$\mathbb{P},H \in X$ such that $X$ thinks that $\mathbb{P} \in L$,
$|\mathbb{P}| < \omega_2$, $H \subset \mathbb{P}$ is $L$-generic and
$\tau \in L[H]$.  Let $\pi:X \rightarrow \bar{X}$ be the
transitivization map. As $otp(X \cap \omega_2) =
(\omega_2)^{\bar{X}}$ is a limit cardinal of $L$, it follows that
$\pi(H) \subset \pi(\mathbb{P})$ is $L$-generic and $\tau \cap X \in
L[\pi(H)]$.  Thus $otp(X \cap \omega_2)$ remains a cardinal in
$L[\tau \cap X]$ as desired.  \end{proof}

\begin{thm}  The following are equiconsistent.
\begin{enumerate}
\item[(a)] $\omega_2$ is inaccessible in $L$ and $\Sigma^{1}_{2}$ are $\omega_1$-Universally Baire.
\item[(b)] There are club many $\sigma \in [\omega_{2}]^{\omega}$ such
that $otp(\sigma)$ is an $L$-cardinal.
\item[(c)] $\omega_2$ is inaccessible in $L$ and there are stationary many 
$\sigma \in [\omega_{2}]^{\omega}$ such
that $otp(\sigma)$ is a cardinal of $L$.
\item[(d)] There is a $\kappa$ which is $L$-large to an
$L$-inaccessible.
\end{enumerate}
\end{thm}
\begin{proof} (a) implies (b) outright by Lemma 2.3.  The argument for (b) implies 
(c) is implicit in the proof of Theorem 2.2.  Assume (c).  Let $g \subset
Col(\omega,<\omega_1)$ be $V$-generic. Let $\kappa =
\omega_{2}^{V}$.  Then in $L[g]$ there is a stationary set of
$\sigma \in [\kappa]^{\omega}$ such that the order type of $\sigma$
is an $L$-cardinal. Thus if $h \subset Col(\omega_{1}, <
\omega_{2}^{V})$ is $L[g]$-generic then in $L[g][h]$ then the
hypotheses of Theorem 3.1 are satisfied so that (a) holds in the
forcing extension described there.  Thus (a), (b) and (c) are equiconsistent. 
Assume (a).  We will show that $\omega_{1}^{V}$ is $L$-large to
$\omega_{2}^{V}$ in $V[g]$ where $g\subset Col(\omega,\omega_1)$ is
$V$-generic.  Let $X \prec H(\theta)$ be countable.  Let $\pi:X
\rightarrow H$ be the transitive collapse.  Let $Y \prec H(\theta)$
with $\pi,H \in Y$ and let $j: Y \rightarrow M$ be its
transitivization.  Note that $j \circ \pi^{-1} = j(\pi)$.  Call this
map $k$.  We have $$k \upharpoonright L_{\omega_{2}^{H}}:
L_{\omega_{2}^{H}} \rightarrow L_{\gamma}$$

\noindent is fully elementary with critical point $\omega_{1}^{H}$
and this map is an element of $M$.  Because $Y$ sees that $H$ is
countable we have

$$k(\omega_{1}) = Y \cap \omega_1
> \omega_{2}^{H}.$$ \noindent Let $g \subset Col(\omega,\omega_1^{H})$ be
$H$-generic.  Let $\alpha < \omega_2^{H}$ be arbitrary and let $x
\in H[g]$ be a real coding a well-ordering of length $\alpha$.  The
sentence asserting the existence of a transitive model of a
sufficient fragment of set theory containing $y$ which sees an
embedding $k:L_{\alpha} \rightarrow L_{\beta}$ with critical point
$\omega_{1}^{H}$ such that $j(\omega_{1}^{H}) > \alpha$ is
$\Sigma^1_2$ in the parameter $x$.  Hence $H[g]$ sees such an
embedding. As $\alpha$ is arbitrary we conclude that $H[g]$ thinks
that $\omega_{1}^{H}$ is $L$-large to $\omega_{2}^{H}$.  Now apply
$\pi$.  To connect (d) back to (a) assume that $\kappa$ is $L$-large 
to some $L$-inaccessible $\lambda$.  Let $g_{\lambda} \subset Col(\omega,<\lambda)$ be 
$V$-generic.  Then $g_{\lambda}$ is also $L$-generic for the same forcing.  By folding 
the embeddings witnessing our hypothesis (d) into countable submodels and collapsing, we 
see that $\kappa$ is $L$-large to $\lambda = \omega_{1}^{L[g_{\lambda}]}$ in 
$L[g_{\lambda}]$.  For ordinals $\gamma < \lambda$ let $g_{\gamma}$ denote 
$g \cap Col(\omega,< \gamma)$.  We claim that in $L[g_{\kappa}]$ there are stationary 
many $\sigma \in [\lambda]^{\omega}$ such that $otp(\sigma)$ is an $L$-cardinal.  It will 
then follow that (c) holds after forcing with $Col(\omega_1,<\lambda)$.  Let 
$f:\lambda^{<\omega} \rightarrow \lambda$ belong to $L[g_{\kappa}]$ be arbitrary.  Let 
$\delta < \lambda$ be a cardinal of $L[g_{\kappa}]$ such that $f[\delta^{<\omega}] \subseteq \delta$.  Let 
$\alpha$ be a regular cardinal of $L[g_{\kappa}]$ below $\lambda$ which is greater than 
$(\delta^{+})^{L}$ so that $f \upharpoonright \delta^{<\omega} \in L_{\alpha}[g_{\kappa}]$.  
In $L[g_{\lambda}]$ there is an elementary embedding 
$$j:L_{\alpha} \rightarrow L_{\beta}$$ \noindent 
with critical point $\kappa$ such that $j(\kappa) > \alpha$.  By standard arguments (using the fact 
that $Col(\omega,<\kappa)$ is $\kappa$-cc in $L$) this embedding extends to a fully elementary 
$$j:L_{\alpha}[g_{\kappa}] \rightarrow L_{\beta}[g_{j(\kappa)}].$$  

\noindent The embedding is defined by 
$j(\mbox{val}(\tau,g_{\kappa})) = \mbox{val}(j(\tau),g_{j(\kappa)})$ and since it extends $j$ we 
will also denote it by $j$.  Now, let $\sigma$ denote the set $j[\delta]$.  Let $z$ be a real in
 $L_{\beta}[g_{j(\kappa)}]$ coding a well-ordering of length $\delta$.  The structure 
$L_{\beta}[\hat{g}]$ has a tree $T$ consisting of pairs $(s,t)$
 with $s$ a finite approximation to a set of ordinals $\sigma$ closed under
 $f$, and $t$ a finite approximation to an order isomorphism between $\delta$
and $\sigma$.  $L_{\beta}[g_{j(\kappa)}]$ must see a branch through this
tree and the result follows by reflection using $j$.   \end{proof}

\begin{lem}  Suppose $\kappa$ is a weakly compact cardinal and $\mathbb{P}$ is
$\kappa$-cc.  Suppose $\Pi^1_2$ sets are $<\kappa$-Universally Baire in $V[G]$
where $G \subset \mathbb{P}$ is $V$-generic. Then $\Pi^1_2$ sets are $\kappa$-cc-UB in $V[G]$.
\end{lem}

\begin{proof}  We will use the fact that a set $A \subset \omega^{\omega}$ is $\kappa$-cc Universally 
Baire iff there are trees $S,T$ on some $\omega \times \lambda$ which project to $A$ and its 
complement and continue to project to complements after forcing with any $\kappa$-cc poset.  So let 
$\dot{\mathbb{Q}}$ be a $\mathbb{P}$-name for a poset which is forced by $\mathbb{P}$ to
have the $\kappa$-cc.  Thus $\mathbb{P} * \dot{\mathbb{Q}}$ has the
$\kappa$-cc in $V$.  Now suppose $\dot{x}$ is a
$\mathbb{P}* \dot{\mathbb{Q}}$-name for a real.  Since $\mathbb{P}*
\dot{\mathbb{Q}}$ is $\kappa$-cc and $\kappa$ is weakly compact
there is an elementary suborder $\mathbb{A} \prec \mathbb{P}*\dot{\mathbb{Q}}$ which has
size strictly less than $\kappa$, decides $\dot{x}$, and has the property that maximal
$\mathbb{A}$ antichains are maximal antichains in $\mathbb{P}*
\dot{\mathbb{Q}}$. The upshot of this is that over $V[G]$ where $G
\subset \mathbb{P}$ is $V$-generic, every real which is generic for
a $\kappa$-cc forcing is generic for a forcing of size $<\kappa$.
Let $A$ be a $\Sigma^1_2$ set and for each forcing $\mathbb{Q}$ of size
$<\kappa$ (whose underlying set is some ordinal below $\kappa$ say)
let $S_{\mathbb{Q}}, T_{\mathbb{Q}}$ be $\mathbb{Q}$-UB
representation of $A$.  These trees may be joined to produce the desired $\kappa$-cc Universal 
Baire representation of $A$.  \end{proof}

\begin{thm} The following are equiconsistent.
\begin{enumerate}

\item $\Sigma^1_2$ sets are $(2^{\omega})^{+}$-cc-Universally Baire

\item $\omega_2$ is weakly compact in $L$ and there are stationary many
$\sigma \in [\omega_2]^{\omega}$ such that $otp(\sigma)$ is a cardinal of $L$

\item There is a $\kappa$ which is $L$-large to an weakly compact cardinal of $L$.
compact.
\end{enumerate}
\end{thm}

\begin{proof}  This is indentical to the proof of Theorem 3.2, using
lemma 3.4 in the argument from (2) to (1) to get the stronger
conclusion.  Of course we are using that $CH$ holds in all models under consideration.  We need to show 
that (1) implies that $\omega_2$ is weakly compact in 
$L$.  Assume (a) and let $\mathbb{Q}$ be the poset for forcing Martin's Axiom.  Let 
$\mathbb{P} = Col(\omega,\omega_1) * \mathbb{Q}$ and note that $\mathbb{P}$ is 
$\omega_2$-cc.  In the extension $V[G]$ by $\mathbb{P}$ we will have $\Sigma^1_2$ sets 
ccc-Universally Baire.  Thus $\Sigma^1_2$ sets are Lebesgue measurable and have the property of 
Baire which implies that $\omega_1$ is inaccessible to reals (see [4]).  Thus in 
$\omega_1 = \omega_{2}^{V}$ is weakly compact in $L$ in $V[G]$ by a result of Harrington and 
Shelah (see [2] or Lemma 7 of [3]).  \end{proof}

\section{Weak Capturing does not imply capturing}

\noindent A weakening of the Universal Baire property is presented in $[6]$. A
set of reals $A$ is weakly captured at $\kappa$ if there is a
$Col(\omega,\kappa)$-term $\dot{A}$ such that for sufficiently large
$\theta$, for a club of countable $H \prec H(\theta)$, and for a
comeager set of $g:\omega \rightarrow otp(H \cap \kappa)$,
$$\pi_{H}(\dot{A})_{g} = A \cap H[g],$$ where $otp(H \cap \kappa)$ is the
order type of $H \cap \kappa$ and $\pi_{H}$ is the transitivization
map.  A less metamathematical characterization is the following.

\begin{lem}  The following are equivalent.

\begin{enumerate}

\item $A$ is weakly captured at $\kappa$.

\item For every continuous $f:\kappa^{\omega} \rightarrow \omega^{\omega}$ there is a dense set of 
$p \in \kappa^{<\omega}$ such that either $f^{-1}(A) \cap \sigma^{\omega}$ is meager below $p$ in 
$\sigma^{\omega}$ for a club of $\sigma \in [\kappa]^{\omega}$ or comeager below $p$ in 
$\sigma^{\omega}$ for a club of $\sigma \in [\kappa]^{\omega}$.
\end{enumerate}
\end{lem}

\begin{proof}  (1) implies (2) is immediate as any condition $p \in \kappa^{<\omega}$ has a refinement 
$p \subseteq q$ such that $q \Vdash_{Col(\omega,\kappa)} f(\dot{g}) \in \dot{A}$ or 
$q \Vdash_{Col(\omega,\kappa)} f(\dot{g}) \notin \dot{A}$.  For the other direction, if $\tau$ is 
a standard $Col(\omega,\kappa)$ term for a real then $\tau$ gives rise to a function 
$f_{\tau}:\kappa^{\omega} \rightarrow \omega^{\omega}$ defined by $f_{\tau}(g) = \tau_{g}$ which 
is continuous on a comeager set.  Define $\dot{A}$ to be 
the set of $(p,\tau)$ such that $f^{-1}(A) \cap \sigma^{\omega}$ is comeager below $p$ in 
$\sigma^{\omega}$ for a club of $\sigma \in [\kappa]^{\omega}$.  A straightforward argument shows that 
$\dot{A}$ is a weak capturing term for $A$. \end{proof}

\begin{thm}  It is consistent that $\Sigma^1_2$ is weakly captured at
$\omega_2$ but not fully captured at $\omega_2$.
\end{thm}

\begin{proof}  Suppose $\kappa$ is $L$-large to an $L$-weakly compact, and
let $L[g][h]$ be the model of Theorem 2.2 in which $\kappa = \omega_1$ and
$\lambda = \omega_2$.  We have shown that there is a stationary set
$S \subset [\omega_2]^{\omega}$ such that whenever $X \prec
H(\theta)$ is such that $X \cap \omega_2 \in S$ then $X[g]$ is
$\Sigma^1_2$ elementary in $V$ for every $X$-generic $g \subset
Col(\omega,\omega_{1} \cap X)$.  The forcing $\mathbb{Q}$ of Lemma 3.1 puts a club 
through $S$ and so in the extension all $\Sigma^1_2$ sets are $\omega_1$-Universally Baire. 
We first argue that $WRP_{(2)}(\omega_2)$ holds in $L[g][h][G]$.  For $a \subseteq \lambda$ let 
$\mathbb{Q}_{a}$ denote the countable support product of $\mathbb{P}_{\gamma}$ taken over 
$L$-cardinals $\gamma \in a$, with $S$ the underlying stationary set.  Returning to $L[g]$ where 
$\lambda$ is still weakly compact, let $p$ be a condition and $\dot{T}$ a term such that 
$$p \Vdash^{L[g]}_{Col(\omega_1,<\lambda)*\mathbb{Q}_{\lambda}} \dot{T}  
\mbox{ is stationary in} [\omega_2]^{\omega}.$$
\noindent By the usual reflection argument we have an inaccessible $\delta < \omega_2$ such that 

$$p \Vdash^{L[g]}_{Col(\omega_1,<\delta) * \mathbb{Q}_{\delta}} 
\dot{T}_{\delta}  \mbox{ is stationary in } [\delta]^{\omega},$$

\noindent where $\dot{T}_{\delta}$ denotes $\dot{T} \cap V_{\delta}$.  Now let 
$h \subset Col(\omega_1,<\lambda)$ be $L[g]$-generic and 
$G \subset \mathbb{Q}$ be $L[g][h]$ generic below the condition $p$.  Let $h_{\delta}$ 
and $G_{\delta}$ be the restrictions to $Col(\omega_1,<\delta)$ and $\mathbb{Q}_{\delta}$ 
respectively.  These are $L[g]$-generic as well and 
$$val(\dot{T}_{\delta},h_{\delta}*G_{\delta}) = val(\dot{T},h*G) \cap [\delta]^{\omega}$$ 
\noindent in $L[g][h][G]$, where.  The key point is that the stationarity of 
$T_{\delta} = int(\dot{T}_{\delta},h_{\delta}*G_{\delta})$ is preserved.  It suffices 
to show that the stationarity of $T$ is preserved by $\mathbb{Q}_{\lambda \setminus \delta}$ 
over $L[h][G_{\delta}]$.  The key point is that 
$\{ \sigma \cap \delta \  | \ \sigma \in S \}$ contains a club in $[\delta]^{\omega}$.  Thus 
if $\dot{C}$ is a name for a club subset of $[\delta]^{\omega}$, we can find a dense set of conditions 
$t \in \mathbb{Q}_{\lambda \setminus \delta}$ with a corresponding $\sigma \in S$ such that 
$\sigma \cap \delta \in T$ and $t$ forces that $\sigma \in \dot{C}$.  Now let $G$ be 
$L[g][h]$-generic for the forcing
$\mathbb{Q}$.  This forcing does not add countable sets of ordinals.
Let $A$ be a $\Sigma^1_2$ set in $L[g][h][G]$.  Then $A =
A^{L[g][h]}$.  Fix such an $A$.  We claim that $A$ is weakly captured in $L[g][h][G]$.   
Otherwise 
there is a condition $t \in \mathbb{Q}$, terms $\dot{T}_{m}$ and $\dot{T}_{c}$, and 
$p \in \omega_{2}^{<\omega}$ such that $t$ forces 

\begin{enumerate}

\item $\dot{T}_{m}$ and $\dot{T}_{c}$ are both stationary subset of $S$

\item $\sigma \in \dot{T}_{m}$ implies $\dot{f}^{-1}(A) \cap \sigma^{\omega}$ is meager below $p$

\item $\sigma \in \dot{T}_{c}$ implies $\dot{f}^{-1}(A) \cap \sigma^{\omega}$ is comeager below $p$.

\end{enumerate}

\noindent We may assume that there is $\delta < \omega_2$ such that 
$t$ forces that both $\dot{T}_{m}$ and $\dot{T}_{c}$ reflect to 
$\delta$.  Let $\bar{t} \leq t$ and $\bar{p} \leq p$ such that 
$$\bar{t} \Vdash_{\mathbb{Q}} \dot{f}^{-1}(A) \cap \delta^{\omega} 
\mbox{ is comeager below } \bar{p}.$$   \noindent It follows that $\bar{t}$ forces that 
$\dot{f}^{-1}(A) \cap \sigma^{\omega}$ is comeager below $\bar{p}$ for a club of 
$\sigma \in [\delta]^{\omega}$, a contradiction.  To finish the proof of the theorem we must 
show that $\Sigma^1_2$ sets are not $\omega_2$-Universally Baire in
$L[g][h][G]$.  Let $D$ be the set of $\alpha < \lambda$ such that
$cf(\alpha) = \omega$ in $L$.  As $L[g][h][G]$ is a $\lambda$-cc
extension of $L$ we know that $D$ remains stationary in $L[g][h][G]$.  Thus the
set of regular cardinals of $L$ below $\omega_2$ cannot be club, as they be would be if 
$\Sigma^1_2$ sets were $\omega_2$-Universally Baire by the argument of Theorem 2.2   \end{proof}

\end{document}